\documentclass[11pt]{article}

\usepackage{latexsym, epsfig}
\usepackage{color}
\usepackage{amssymb, amsmath, amsfonts}
\usepackage{exscale}
\usepackage{ifthen}
\usepackage{longtable}
\usepackage{subfigure}
\usepackage{caption}
\usepackage{cite}
\usepackage{lscape}
\usepackage{enumitem}

\newcommand{\vol}{\mbox{\rm vol}}

\newcommand{\qed}{$\hfill\Box$}

\newtheorem{definition}{\bf Definition}[section]

\newtheorem{theorem}[definition]{\bf Theorem}
\newtheorem{lemma}[definition]{\bf Lemma}
\newtheorem{proposition}[definition]{\bf Proposition}

\newtheorem{remark}[definition]{\sc Remark}

\begin{document}
\bibliographystyle{alphabetical}

\protect\pagenumbering{arabic}

\title{\sc Entropy Rigidity of negatively curved manifolds of finite volume}
\author{\sc M. Peign\'e, A. Sambusetti}
\date{\today}
\maketitle

\centerline {\small {\bf Abstract}}
We prove the following entropy-rigidity result in finite volume: if $X$ is a negatively curved manifold with curvature $-b^2\leq K_X \leq -1$, then   $Ent_{top}(X) = n-1$ if and only if $X$ is hyperbolic. In particular, if $X$ has the same length spectrum of a hyperbolic manifold $X_0$, the it is isometric to $X_0$  (we also give a direct, entropy-free proof of this fact).  We  compare with the classical theorems holding in the compact case, pointing out the main difficulties to extend them to finite volume manifolds. 
\bigskip


 \section{Introduction}
 
 The problem of length spectrum  rigidity of Riemannian manifolds has a long history. 
 The fact that, in negative curvature  (even in constant curvature), the collection of the lengths of all closed geodesics, together with all multiplicites, does not determine  the metric is well known since \cite{vigneras}. 
 On the other hand, on a {\em compact}, negatively curved surface $\bar X$, the metric is determined up to isometry by the {\em marked length spectrum} (that is, the map ${\cal L}: {\cal C}(\bar X) \rightarrow \mathbb{R}$ associating to each free homotopy class of loops in $\bar X$ the length of the shortest geodesic in the class); this was proved by Otal \cite{otal} and, independently, by Croke \cite{croke}. 
 The same is true in dimension $n \geq 3$ for any compact, {\em locally symmetric} manifold $\bar X_0$ of negative curvature: the locally symmetric metric on $\bar X_0$ is determined, among all negatively curved metrics, by its marked length spectrum. This is  consequence of  Besson-Courtois-Gallot's solution of the minimal entropy conjecture  and of the fact, proved by Hamenstadt \cite{hamen1}, that if a compact, negatively curved manifold $\bar X$ has the same marked length spectrum as a compact,  locally symmetric space $\bar X_0$, then ${\rm vol}(\bar X)={\rm vol}(\bar X_0)$
 \footnote{  By \cite{hamen2}, two compact, negatively curved manifolds having the same marked length spectrum  have $C^0$-conjugated geodesic flow; moreover, if a compact manifold $\bar X$ has geodesic flow which is $C^0$-conjugated to the flow of a manifold $\bar X_0$ whose unitary tangent bundle  has a $C^1$-Anosov splitting (e.g., a locally symmetric space), then $\bar X$ has the same volume as $\bar X_0$,  see \cite{hamen1}. The fact that, for compact manifolds, volume is preserved under $C^1$-conjugacies is much easier and relies on Stokes' formula, cp. \cite{ck}.}.
 
Less seems known about the length rigidity of negatively curved,  {\em finite-volume} manifolds: most generalizations are not straightforward, and seem to require  additional  assumptions (such as bounds on the curvature and  on its derivatives, or  the finiteness of the Bowen-Margulis measure); we will try to point out some of these difficulties throughout the paper. 
 For instance, the fact that having the same marked length spectrum  implies the existence of a $C^0$-conjugacy of the geodesic flow, would certainly require  some new arguments  for finite-volume manifolds 
\footnote{Cp. Lemma 2.4 in \cite{hamen2}, which is central in the argument: it is based on the fact that   the closed geodesics on $\bar X$ equidistribute  towards the Bowen-Margulis measure when $\bar X$ is compact. This property does not   hold for non uniform lattices  as soon as the Bowen-Margulis measure is infinite.
 }.  
 \vspace{-2mm}

 \begin{theorem}
 \label{teor1}
 Let $\bar X$  be a finite volume $n$-manifold with pinched,  negative curvature $-b^2\leq K_{\bar X} \leq -1$ which is homotopically equivalent to a locally symmetric manifold $\bar X_0$,\linebreak with curvature normalized between $-4$ and $-1$.   If $\bar X$ and $\bar X_0$ have same marked length spectrum, then they are isometric.
 \end{theorem}

\vspace{-2mm}
The proof of this is probably known to experts and follows a classical scheme: one can construct a $\Gamma$-equivariant map $f: X \rightarrow X_0$ between the universal coverings, which induces a homeomorphism between the boundaries and preserves the cross-ratio; then, the conclusion stems, for instance, from Bourdon's result  \cite{bourdon2} on M\"obius embeddings from  locally symmetric to $CAT(-1)$-spaces.
However, the main difficulty, in the case of finite volume manifolds,  is to show that $f$ is a quasi-isometry, the quotients   $\bar X$ and $\bar X_0$ being non-compact; we will give a short proof of this fact in \S\ref{sectiontheorem1},  by way of example, to measure the difference from the compact case. 

 
  It is tempting to approach the above problem  by using a finite-volume version of Besson-Courtois-Gallot's inequality, given by Storm \cite{storm}; however, notice that Storm's inequality  
$Ent(\bar X)^nVol(\bar X) \geq Ent(\bar X_0)^n Vol(\bar X_0)$ 
concerns the {\em volume entropy}\footnote{Cp. the definition of the maps $\Psi^b_c$ in \cite{storm}, which clearly require that $c$ is greater than the exponential growth rate of the universal covering of the manifold under consideration.} of $\bar X$, and not  the {\em topological entropy}  $Ent_{top}(\bar X)$ of the geodesic flow on $U\bar X$.
Recall that  for compact, negatively curved manifolds,  one always has  $Ent(X) =Ent_{top} ( \bar X) $, but for finite-volume manifolds   $Ent(\bar X)$ is generally strictly greater than $Ent_{top} ( \bar X) $ (cp. \cite{DPPS-crelle}, \cite{DPPS-fourier}); on the other hand, the topological entropy  always equals the   {\em critical exponent} of the group
$\Gamma = \pi_1(X)$ acting on the universal covering:
  \vspace{-3mm}
  
  $$\delta_{\Gamma}:= \lim_{R \rightarrow \infty}  {1\over R}  \ln \#  \{ \gamma \in \Gamma \; | \; d(x,  \gamma   x) \leq R \} $$

  \vspace{-1mm}
\noindent as proved in \cite{otalpeigne}.
Then, the volume entropy $Ent(\bar X)$ is not   preserved,   a priori,  by the condition of having same marked length spectrum, or by a conjugacy of the geodesic flows. Moreover, it is not  clear whether, for finite volume manifolds, the volume is preserved under a conjugacy  of the flows
\footnote{This seems unclear even under the assumption of a $C^1$-conjugacy;  cp. the proof of Proposition 1.2 in \cite{ck}, where Stokes's theorem fails, unless one knows that the conjugacy $F$ has bounded derivatives.}.
 
The  upper bound on the curvature  $K_{\bar X} \leq -1$ in Theorem \ref{teor1} seems unreasonably strong,  as it implies,  when $\bar X_0$ is hyperbolic, that  $\bar X$  has  marked length spectrum which is {\em asymptotically critical}: that is, its exponential growth rate $\delta_\Gamma$ is greater than or equal to the corresponding exponential growth rate for $\bar X_0$ (cp. \cite{DPPS-crelle}, Lemma 4.1). \linebreak
We expect that the same result holds without curvature bounds,   but, even in this weaker form,  we were unable to find a proof of this  result in literature.

The knowledge of the full marked length spectrum can be relaxed, as we show in the following result (which implies Theorem \ref{teor1} in the real hyperbolic case):


\begin{theorem}
\label{teor2} 
Let $\bar X$  be a finite volume $n$-manifold with pinched,  negative curvature $-b^2\leq K_{\bar X} \leq -1$.
Then  $Ent_{top}( \bar X) \geq  n-1$, and the equality $Ent_{top}( \bar X) = n-1$ holds if and only if $ \bar X$ is hyperbolic. 
\end{theorem}

 The  entropy  characterization of constant curvature (and locally symmetric) metrics has been declined in many different ways so far:  in the compact case, the above theorem is due to   Knieper (see   \cite{knieper1}, where this result is not explicitly stated,  but can be established following the argument of the proof of Theorem 5.2.); see also  \cite{courtois}, and  \cite{BK} for a proof in the convex-cocompact case.

 We want to stress here that a basic difference between  Theorem \ref{teor2} 
(or, more precisely,   their compact versions in  \cite{knieper1},\cite{courtois},\cite{BK}) 
and  the celebrated entropy characterization  of  Hamenst\"adt  \cite{hamen0}   of locally symmetric metrics,  with same curvature normalization,   is the lack of any locally  symmetric manifold $\bar X_0$  of reference homotopically equivalent to $\bar X$. 
Actually, the characterization given by  Theorem \ref{teor2} is very particular to  constant curvature spaces and 
it does not generalize, as it is, to  locally symmetric spaces: indeed,  it is easy to construct compact, pinched, negatively curved manifolds with $-b^2 \leq K_X \leq -1$ having same entropy as, let's say, the complex hyperbolic space, but which are not complex hyperbolic.

The same  difference holds with the existing, finite volume  versions of Besson-Courtois-Gallot's characterization of locally symmetric spaces, in particular with \linebreak Boland-\-Connell-Souto's  papers  \cite{BCS} and  Storm's \cite{storm}: these two works,   together,  imply that if a  finite volume  manifold $\bar X$  with  curvature  $K_{\bar X} \leq -1$ has volume entropy $Ent (\bar X)=n-1$,  then it is hyperbolic,  {\em provided that  one knows beforehand  that $\bar X$ is   homotopically equivalent to a  hyperbolic manifold} $\bar X_0$.
Besides the difference between volume and topological entropy stressed above,   this strong  supplementary  topological assumption on   $\bar X$ is not  made  in   Theorem \ref{teor2}.

Let us also  point out that   Knieper's approach in \cite{knieper1} does not allow to deduce the above characterization in the finite volume case. Although  G. Knieper's horospherical measure $\mu_H$ can 
be perfectly defined in this context (following \S3 of  \cite{knieper1}),  it can easily be infinite,  as well as the Bowen-Margulis measure $\mu_{BM}$: given a finite volume surface $\bar X$ with  convergent fundamental group $\Gamma$ and with a cusp whose metric, in horospherical coordinates, writes as ${\cal A}^2(t) dx^2 + dt^2$, it is not difficult to show that $\mu_H$  is infinite as soon as 
$$ \int_0^\infty e^{\delta_\Gamma t} {\cal A}(t)dt =+\infty$$
(cp. examples in \S3, \cite{DPPS-fourier}). Therefore, all formulas in \cite{knieper1} relating $Ent_{top} (\bar X)$ to  the trace of the second fundamental form of unstable horospheres need to be justified in some other way 
\footnote{For instance, Corollary 4.2 in   \cite{knieper1} only holds for $\mu_H$-integrable functions, and cannot be applied as it is to constant functions or to $tr \, U^+(v)$ to deduce Theorem 5.1, when $\| \mu_H\| =\infty$.}.

%
%
%

On the other hand, we will give in \S\ref{Entropyrigidity} a proof of  Theorem \ref{teor2} using the barycenter method,  initiated by  Besson-Courtois-Gallot in \cite{bcg}-\!\!\cite{bcg2}, together with some  careful  estimates of the Patterson-Sullivan measure, which will not need neither the finiteness of $\mu_{BM}$ (or $\mu_H$) nor the conservativity of the geodesic flow with respect to $\mu_{BM}$.


Also,  notice that if we drop the assumption $K_{\bar X} \geq -b^2$ in Theorem \ref{teor2},  the manifold $\bar X$ might as well be of infinite type (i.e. with infinitely generated fundamental group,  or even without any cusp,  see examples in \cite{nguyen}),  hence very far from being  a hyperbolic manifold  of finite-volume.
\vspace{4mm}
 
 \small
 {\sc Notations}. \\
Given  functions $f,g: \mathbb{R}_+ \rightarrow \mathbb{R}_+$, we  will systematically write  
  $f\stackrel{C}{\prec}g $ (or $g \stackrel{C}{\succ} f$)  if   there exists $C>0$ and $R_0>0$ such that  $f(R) \leq C g( R)$  for $R>R_0$.  We  write $f  \stackrel{C}{\asymp} g$  when  $g  \stackrel{C}{\prec}  f \stackrel{C}{\prec}  g$ for $R \gg0$ (or simply   $f \asymp g$  and $f \prec g $   when the constants $C$ and $R_0$ are unessential) 
 \normalsize

 \section{Geometry at infinity in negative curvature}

 Throughout all the paper,  $X$ will be  a  $n$-dimensional, complete,  simply connected manifold with  strictly negative curvature $-b^2\leq K_X \leq -1$.
 
Let  $X(\infty)$  the ideal boundary of $X$: for  $x, y \in X$ and $\xi \in  X (\infty)$,  we will denote by $[x,y]$ (resp. $[x, \xi[$) the geodesic segment  from $x$ to $y$ (resp. the ray  from $x$ to $\xi$), and by $x\xi(t)$ the parametrization of geodesic ray $[x, \xi[$ by arc length.  Let
$$b_{\xi}(x, y)=\lim_{z\to\xi}d(x, z)-d(z, y)$$
 be the Busemann function centered at $\xi$; 
the level set $  \partial H_\xi (x)= \{ y \,  | \,  b_\xi (x, y)\vert= 0\}$ (resp. the suplevel set \nolinebreak $H_\xi (x) = \{ y \,  | \,  b_\xi (x, y)  \geq  0\}$ is the horosphere (resp. the horoball) with center  $\xi$ and passing through $x$. 
We will denote by $d_\xi$  the horospherical distance between two points on a same horosphere centered at $\xi$, and we define the   radial semi-flow  $(\psi_{\xi, t})_{t\geq 0}$  in the direction of $\xi$ as follows: for any $x \in X$,
the point  $\psi_{\xi, t}(x)$  lies on the geodesic ray $[x,\xi[ $  at distance $t$ from $x$. 

Finally, recall that for any fixed $x \in X$,  the Gromov product between  two points  $\xi, \eta \in X(\infty), \xi \neq \eta$, is  defined as 
$$(\xi\vert \eta)_{x } = {  b_{\xi}(x, y)+b_{\eta}(x, y)\over 2}$$
 where $y$ is any point on the geodesic $ ]\xi, \eta[$ joining $\xi$ to $\eta$; as $K_X \leq -1$,  the expression
$ D_x(\xi, \eta)= e^{ - (\xi\vert \eta)_{x} }$
defines  a distance on $X(\infty)$, which we will   call the {\em visual distance  from $x$}, cp. \cite{bourdon1}. 
Accordingly,  the {\em cross-ratio} on $X(\infty)^4$  is defined as 
$$ [\xi_1, \xi_2, \xi_3, \xi_4]={D_x (\xi_1, \xi_3)D_x(\xi_2, \xi_4)\over D_x (\xi_1, \xi_4)D _x(\xi_2, \xi_3)} 
= \hspace{-3mm} \lim_{\stackrel{  p_1, p_2, p_3, p_4 \in X}{  (p_1, p_2, p_3, p_4) \to (\xi_1, \xi_2, \xi_3, \xi_4)}}
\hspace{-5mm} e^{d(p_1,p_3)+d(p_2, p_4)-d(p_1,p_4)-d(p_2,p_3)}$$
for all  $\xi_1, \xi_2, \xi_3, \xi_4 \in X(\infty)$,  and it is easily seen that it is independent from the choice of the base point $x$, cp.\cite{otal-iberico}, \cite{bourdon1}.

  We will  repeatedly make use of the following,  classical result in strictly negative curvature:  
  there exists $\epsilon (   \vartheta)=  \log (\frac{2}{1 - \cos \vartheta})$ such that any geodesic triangle   $xyz$   in $X$ making  angle $\vartheta=\angle_z (x, y) $ at  $z$ satisfies: 
 \vspace{-3mm}
 
\begin{equation}
\label{eqoppositetriangle} d(x, y) \geq d(x, z) + d(z, x) - \epsilon( \vartheta).
\end{equation}

%
%
%

\subsection{ On the geometry of finite volume manifolds}\label{subsectiondecomposition}
 \vspace{-2mm}
Consider  a lattice $\Gamma$ of $X$. The quotient manifold $\bar X= \Gamma \backslash  X $  has finite volume, it is thus a geometrically finite manifold  which admits some particular decomposition which we now recall. The   following classical results are  due to B. Bowditch \cite{bow}, and we state them in the particular case of  finite volume manifolds :
\vspace{1mm}        

\noindent   (a) the limit set of  $L(\Gamma)$ of $\Gamma$ is the full boundary at infinity  $X (\infty)$ and is the disjoint 
union of the radial limit set $L_{rad}(\Gamma)$ with  finitely many orbits
 of {\em bounded} parabolic fixed 
points $ L_{bp} (\Gamma) =  \Gamma \xi_{1} \cup  \ldots \cup \Gamma \xi_{l}$; this means that each  $\xi_i \in L_{bp}(\Gamma) $  is the fixed point of some  maximal   parabolic subgroup  $P_i$ of $\Gamma$, acting   co-compactly on $X(\infty)\small\setminus \{\xi_i\}$;
\vspace{1mm}

\noindent    (b) {\em (Margulis' lemma)}  there exist  closed horoballs $H_{\xi_{1}},  \ldots,  H_{\xi_{l}}$ centered respectively  at $\xi_{1},  \ldots,  \xi_{l}$,  such that  $\gamma   H_{\xi_{i}} \cap H_{\xi_{j}} = \emptyset$  for  all $1\leq i, j\leq l$ and  all $\gamma \in \Gamma \small
 \setminus P_i$;
\vspace{1mm}   
        
\noindent   (c)   The finite volume manifold $\bar X$  can  be decomposed into
 a disjoint union of a compact set $\bar  {\cal K}$ and finitely many ``cusps'' $\bar {\cal C}_1,  ...,  \bar {\cal C}_l$: each $\bar {\cal C}_i$ is isometric to the quotient of $H_{\xi_{i}}$  by a corresponding maximal bounded parabolic group $P_i$. We refer to $\bar  {\cal K}$
and to  $\bar {\cal C} = \cup_i \bar {\cal C}_i$ as to the {\em thick part}  and the {\em cuspidal part} of $\bar X$.
 \vspace{1mm}          

For any fixed $x \in X$, let ${\cal D}={\cal D}(\Gamma, x)$  
the Dirichlet  domain of $\Gamma$  centered at $x $; this is a   convex fundamental subset of $X$,  and we may   assume that  ${\cal D}$ contains the geodesic rays $[x, \xi_i[$.  Each parabolic group $P_i$ acts co-compactly on the horosphere $\partial H_{\xi_i}$ which bounds the horoball $H_{\xi_i}$;  setting 
 ${\cal S}_i = {\cal D} \cap   \partial H_{\xi_i}$ and
 ${\cal C}_i = {\cal D} \cap H_{\xi_i} \simeq   {\cal S}_i \times \mathbb{R}_+$, 
the fundamental  domain ${\cal D}$ 
 can be decomposed into a disjoint union       
 \vspace{-3mm}  
 
\begin{equation}\label{decompositionfundamentaldomain}{\cal D} =  {\cal K}  \cup {\cal C}_1 \cup  \cdots\cup {\cal C}_l
\end{equation}

\vspace{-1mm}  
\noindent where  $ {\cal K}$     is a convex,  relatively compact set containing $x$ in its interior  and projecting to the thick part $\bar  {\cal K}$ of $\bar X$ ,  while   ${\cal C}_i$  and  ${\cal S}_i$  are,  respectively,  connected fundamental domains for the action of $P_i$ on $H_{\xi_{i}}$ and  $  \partial H_{\xi_{i}}$, projecting respectively to    $\bar {\cal C}_i$ and    $\bar {\cal S}_i$.

\vspace{-2mm}  
 \subsection{Growth of parabolic subgroups }
 \vspace{-2mm}
 
   The subgroups $P_1,  s, P_l$ will play a crucial role in the sequel;    
 the growth of their orbital functions is best expressed by introducing the horospherical area function. \linebreak
 Let us recall the necessary  definitions:
 \vspace{-1mm}

\begin{definition}[Horospherical Area]
\label{defiparabolic} ${}$\\
\noindent Let $P$ be a {\em bounded} parabolic group of isometries of $X$ fixing $\xi \in X (\infty)$: that is, $P$ acts cocompactly on $X (\infty) \smallsetminus \{ \xi \}$ (as well as on every horosphere  centered at $\xi$).   \linebreak
Given $x\in X$, let  
  ${\cal S}_x$ be  a  fundamental, relatively compact domain  for the action  of  $P$ on  $  \partial H_\xi (x)$: 
  the  {\em horospherical area  function}  of $P$ is the function 
  \vspace{-3mm}
    
$${\cal A}_P (x, R) = \vol \left[ P  \backslash \psi_{\xi, R} \left(  {  \partial H_{\xi} (x)}  \right) \right]  
= \vol \left[   \psi_{\xi, R} \left(  {\cal S}_x   \right) \right]  $$
 where  $\vol$ denotes the Riemannian measure of horospheres. 
  \end{definition}
 \pagebreak
  
 \begin{remark}
{\em 
 When $-b^2 \leq K_X \leq -a^2 <0$, well-known estimates  of the differential of the radial flow (cp. \cite{HIH}) 
 yield,  for any $t \in \mathbb R $ and $v \in T^1 X$

\begin{equation}
\label{jacobi-norm}
e^{-bt}  \parallel   v    \parallel \leq  \parallel    d\psi_{\xi,  t}(v)     \parallel \leq e^{-at} \parallel     v    \parallel
\end{equation}
Therefore we deduce that,  for any $\Delta> 0$,

\begin{equation}
\label{eqA(R+C)}
e^{-(n-1)b\Delta}\leq  \frac{ {\mathcal A}_P(x,  R+\Delta) }{{\mathcal A}_P(x,  R) } \leq e^{-(n-1)a\Delta}
\end{equation}
}
 \end{remark}

The following Proposition  shows   how the horospherical area  ${\cal A}_P$ is   related to the orbital  function of $P$, cp. \cite{DPPS-crelle}:

\begin{proposition}
\label{propVp}
Let $P$ be a bounded parabolic group of $X$ fixing $\xi$,  with  $diam ({\cal S}_x) \leq d$. \\
There exist    $R_0$ and  $\Delta_0$
only depending on $n, a, b, d$ and constants $C=C(n,  a, b,  d)$ and  $C'=C'(n,  a, b,  d, \Delta)$ such that,  for any $R \geq   b_\xi(x, y)  +R_0$ and any $\Delta> \Delta_0$, the numbers $v_P (x, y,  R)$ and  $v^{\Delta}_P (x, y,  R)$ of orbit points of $Py$ falling, respectively, in the balls $B(x,R)$ and in the annuli $A^\Delta(x,R)$ satisfy:
 \vspace{-3mm}

$$v_P (x, y,  R)   = \{ p \in P \, | \,  d(x,py) < R \}  \stackrel{C}{\asymp}    {\cal A}_P^{-1} \left(x,  \frac{R+ b_\xi (x, y)}{2} \right)$$
 
 \vspace{-6mm}
 
$$v^{\Delta}_P (x, y,  R) = \{ p \in P \, | \,  R - \frac{\Delta}{2} \leq d(x,py) \leq  R + \frac{\Delta}{2} \}    \; \stackrel{C'}{\asymp}  \;  {\cal A}_P^{-1} \left(x,  \frac{R+ b_\xi (x, y)}{2} \right).$$

\end{proposition}

 \section{Length spectrum and rigidity}\label{Lengthspectrumandrigidity}
\label{sectiontheorem1}
This section is devoted to the proof of Theorem \ref{teor1}. \\
 Let $\Gamma$ be the fundamental group of the manifolds $\bar X$ and $\bar X_0$, acting by isometries on their Riemannian universal coverings $X$ and $X_0$ respectively. We will construct a $\Gamma$-equivariant homeomorphism 
 $f_{\infty}: \partial X(\infty) \rightarrow X_0 (\infty)$ and apply the following: 
 
 \begin{theorem}[cp.  \cite{bourdon2}]
 \label{bourdon}
 Let $X$ be a {\rm CAT}$(-1)$-space and $X_0$ a symmetric  space of rank one, with curvature $-4 \leq K_{X_0} \leq 1$.
Assume that  $f_{\infty}: X(\infty) \rightarrow X_0(\infty)$  is  a $\Gamma$-equivariant homeo\-morphism which preserves the cross-ratio: then there exists a $\Gamma$-equivariant  isometry $f: X \to X_0$  whose  extension on $X(\infty)$ coincides with $f_{\infty}$.
 \end{theorem}
 
\noindent For this, we fix $x \in X$ and  $x_0 \in X_0$ and consider the natural  $\Gamma$-equivariant bijection  $\phi: \Gamma  x \rightarrow \Gamma   x_0$. The main difficulty here is to show the following: 
 
 \begin{proposition}
 \label{prop-quasi}
 The map $\phi$ is a quasi-isometry between the orbits, with respect to the distances induced by the Riemannian distances of $X$  and $X_0$ respectively.
 \end{proposition}
 
 We assume Proposition \ref{prop-quasi} for a moment.  Since $\bar X$ and $\bar X_0$ have finite volume, the  limit set of $\Gamma$ coincides with the full boundaries  $X(\infty)$ and $X_0(\infty)=\mathbb S^{n-1}$ and the map $\phi$ extends to a  bi-H\"older  and $\Gamma$-equivariant homeomorphism $f_\infty$  between these boundaries, endowed with their natural visual metric from $x$ and $x_0$.   \\
  Now,  the fact that $\bar X$ and $\bar X_0$ have the same marked length spectrum implies that $f_\infty$ preserves the cross ratio; this follows  for instance from \cite{otal-iberico}.  For the sake of completeness, we will give a  proof  of this fact at the end of this section (Proposition \ref{crossratio}), based on an argument from \cite{kim} (where the same is proved for symmetric spaces).\\
 We  conclude by  \ref{bourdon} that there exists an isometry between the quotients   $\bar X$ and $\bar X_0$.\qed
\vspace{5mm}

 {\bf Proof of Proposition \ref{prop-quasi}.}\\
Let us first show that there exists $\lambda>1$ such that, for all $\gamma \in \Gamma$, we have
\begin{equation}\label{quasiiso-1sense}
d_0(x_0, \gamma  x_0)\leq \lambda d(x, \gamma   x) +\lambda.
\end{equation}
Consider the decomposition of  $\bar X$ described in subsection  \ref{subsectiondecomposition}: 
we denote by  $ \mathcal H $ the set  of pairwise disjoint horoballs which project on the cuspidal part of $\bar X$, so that $ \widetilde{ \mathcal K} := X\small\setminus \cup_{H\in \mathcal H} H =  \Gamma {\mathcal K} $ is the subset of $X$ projecting to the thick part   $\bar {\mathcal K}$ of $\bar X$.\\ 
We assume that    $d( H, H') \geq 1$ for any $H\neq H'$ in $\mathcal H$, and set ${\rm diam}({\cal K})= D$.  \\
For any $\gamma \in \Gamma$, the geodesic  segment   $[x, \gamma  x]$  intersects at most $N \leq d(x, \gamma   x) $  distinct horoballs $ H\in \mathcal H$ and can be decomposed as 
\vspace{-3mm}

$$
[x, \gamma   x]= [x_0^+, x_1^-]\cup[x_1^-, x_1^+]\cup   \cdots\cup[x_{N-1}^+, \gamma  x_N^-]
$$

\vspace{-1mm}
\noindent  where $x_0^+=x$, $x_N^-=x$, and where  $[x_i^-, x_i^+]$ is equal to $[x, \gamma   x]\cap H_i$ for some horoball $H_i\in \mathcal H$ and each $[x_i^+, x_{i+1}^-]$ is included in $\widetilde{ \mathcal K}$.
Then, there exist  elements $g_i \in \Gamma$ and $p_i \in P_1 \cup  \cdots\cup P_l$, for $ 1 \leq i\leq N-1$,  with $g_N := \gamma$, such that
$x_i^-\in g_i   {\cal K}$, $x_i^+\in g_i p_i   {\cal K}$; moreover, set $\gamma_i:= p_{i-1}^{-1} g_{i-1}^{-1}g_i$ for $1\leq i \leq N$ with the convention $p_0=g_0=1$. \\
Notice that all the geodesics  $[  x ,   \gamma_{i} x ]$ are included in a $D'=D'(D)$-neighbourhood of $\widetilde {\cal K}$: \linebreak
in fact, the length of the broken geodesic   $[x_i^+,  g_i p_i x ] \cup [ g_i p_ix ,   g_{i+1}x ] \cup [ g_{i+1}x,  x_{i+1}^-] $ exceeds the length of
$[x_i^+,  x_{i+1}^-] $ at most of $2D$, so (the curvature being bounded above by $-1$) it stays $D'(D)$ close to $[x_i^+,  x_{i+1}^-] $;  by construction this last segment does not  enter any  horoball  of $ {\cal H}$, so $[g_{i}p_{i} x, g_{i+1} x]$ and   $ [  x ,   \gamma_{i+1} x ] =   p_i^{-1}g_{i}^{-1} [g_{i}p_{i} x, g_{i+1} x]$    stay   in the $D'$-neighbourhood of $\widetilde {\cal K}$.

\noindent Now we have $\gamma= \gamma_1 p_1\gamma_2  \cdots\gamma_{N-1}p_{N-1}\gamma_N$, 
so
\vspace{-3mm}

\begin{equation}\label{ineqtriangu1}
 d_0(x_0, \gamma  x_0) \leq \sum_{i=1}^N d_0(x_0, \gamma_i x_0)+\sum_{i=1}^{N-1} d_0(x_0, p_i  x_0).
\end{equation} 

\vspace{-3mm}
\noindent On the other hand 
\vspace{-6mm}

\begin{eqnarray*} 
 d(x, \gamma  x)
 &=&   \sum_{i=1}^{N}d(x_{i-1}^+, x_{i}^-)  +  \sum_{i=1}^{N-1} d(x_i^-, x_i^+)
  \\
   &\geq&   \sum_{i=1}^{N-1} d(g_{i-1}p_{i-1} x, g_{i}  x) + \sum_{i=1}^{N-1} d(g_i  x, g_ip_i x) -4(N-1)D\\
   &=& \sum_{i=1}^N d(x, \gamma_i x) + \sum_{i=1}^N d(x, p_i  x)-4(N-1)D\\
\end{eqnarray*} 

\vspace{-3mm}
\noindent which in turn yields, as  $N  \leq d(x, \gamma  x)$,  
\begin{equation}\label{ineqtriangu2}
 \sum_{i=1}^N d(x, \gamma_i  x)+\sum_{i=1}^N d(x, p_i  x)\leq  (1+4D)d(x, \gamma  x)  .
\end{equation}
To obtain  inequality (\ref{quasiiso-1sense}), it is thus sufficient to check that it holds for each $\gamma_i$ and $p_i$  which appears in the sums (\ref{ineqtriangu1}) and  (\ref{ineqtriangu2}).  This is proved in  the  two following lemmas: 

\begin{lemma}\label{quasiisothickpart}
For any $D'>0$, there exists  $C>0$ such that 
\vspace{-3mm}

$$d_0(x_0, \gamma  x_0)\leq C \ d(x, \gamma   x)+C$$

\noindent for any $\gamma \in \Gamma$ such that 
$[x, \gamma   x]$ lies in the $D'$-neighbourhood of   $\widetilde {\cal K}$.
\end{lemma}

\begin{lemma}\label{quasiisohorosphere} There exists   $C'>0$ such that,  for any parabolic isometry $p \in P_1\cup  \cdots \cup P_l$ 
\vspace{-7mm}

\begin{equation}
\label{eq-quasi}
d(x, p  x)\stackrel{C'}{\asymp} d_0(x_0, p  x_0).
\end{equation}
\end{lemma}

\vspace{2mm}
 Switching the roles of $(X, d)$ and $(X_0, d_0)$,  we obtain  the opposite inequality \linebreak $d_0(x_0, \gamma  x_0)\leq \lambda d(x, \gamma   x) +\lambda$, which concludes the proof of  Proposition \ref{prop-quasi}.\qed

\vspace{5mm}
{\bf Proof of Lemma \ref{quasiisothickpart}.} \\
Let  $\gamma \in \Gamma$ such that 
$[x, \gamma   x]$ lies in the $D'$-neighbourhood of   $\widetilde {\cal K}$,   and let $x_0=x$, $x_N= \gamma x$ and  $ x_1, \cdots, x_{N-1}$ be the points on the geodesic segment  $[x, \gamma   x]$ such that $d(x, x_i)= iD$ for $0\leq i \leq N-1$, 
with $N-1= [d(x, \gamma   x)]$.
There exist isometries  $h_0=1,h_1, \cdots, h_{N-1}, h_N=\gamma$ in $\Gamma$ such that 
$d(x_i, h_i   x)\leq D +D'$;  setting   $k_{i} = h_{i-1}^{-1}h_i$,  we then have  $  \gamma=k_1 k_2 \cdots k_{N}. $  
Now, for any $1\leq i \leq N$, we have   $d(x, k_i  x) \leq 1+D+D'$; so every $k_i$ belongs to the finite set $B:=   \{k \in \Gamma \, |\,  d(x, k  x) \leq  1+D+D' \} $.  Setting $C:= \max\{d_0(x_0, k  x_0) \, | \,  k \in B \}$, we obtain
\vspace{-3mm}

$$
d_0(x_0, \gamma   x_0) \leq \sum_{i=1}^{N} d_0(x_0, k_i  x_0)   \leq     N C
   \leq  C d(x, \gamma   x)  +C .\;\;\Box
$$

{\bf Proof of Lemma \ref{quasiisohorosphere}.} \\
Let us first notice that if $p\in \Gamma$ acts on   $X$  as a  parabolic (resp. a hyperbolic) isometry, then it  acts  in the same way on $X_0$: actually,   the infimum of the length of curves in $\bar X$  in the free homotopy class  of a parabolic element $p$  is  $0$ and this condition is preserved  since $\bar X$ and $\bar X_0$ have the same length spectrum.\\
	Then, let $\xi_1, ..., \xi_l \in X(\infty)$ be the fixed points of the maximal parabolic subgroups $P_1,..,P_l$ of $\Gamma$ such that  the geodesic rays $[x,\xi_i[$ are included in the Dirichlet domain ${\cal D}$, as  described  in the subsection \S\ref{subsectiondecomposition}, and call $\xi'_i$ the corresponding parabolic fixed points of $ X_0 (\infty)$; in order to simplify the notations, we set $P=P_i$, $\xi=\xi_i$ and $\xi'=\xi_i'$. \\
Fix a finite generating set $S $ for $P $ and let  $\vert \cdot \vert_{S }$ be the corresponding word metric. \linebreak
As $P $ acts cocompactly by isometries on $(\partial H_{\xi} (x), d_\xi)$ and on $(\partial H_{\xi'} (x_0), d_{\xi'})$ 
we know that these metric spaces are both quasi-isometric to $(P, \vert \cdot  \vert _{S})$. In particular,  there exists a constant $c>0$ such that, for any $ p \in P$
\vspace{-3mm}

\begin{equation} \label{quasiisohorospheres}
{1\over c} d_{ \xi}(x, p  x) -c \leq
 d_{\xi'}(x_0, p  x_0)
\leq
c d_{\xi}(x, p  x) +c.
\end{equation}

\vspace{-1mm}
\noindent Now, by the bounds on curvature $ -b^2 \leq K_X \leq -1$ we get (cp. \cite{HIH})
\vspace{-3mm}

$$  2 \sinh \left( \frac{d (x,px)}{2} \right)    \leq d_\xi (x,px) \leq  \frac{2}{b} \sinh \left( \frac{b}{2} d  (x,px)  \right) $$

\vspace{-2mm}
\noindent hence $ d(x,px) / d(x_0,px_0)  \stackrel{C'}{\asymp} 1$ 
for a constant $C'>0$ only depending on $b$ and $c$.$\Box$

%

\begin{proposition}\label{crossratio} 
Let  $\alpha$ and $\beta$ be two hyperbolic isometries in $\Gamma$ with, respectively,   repelling and attractive fixed points $\alpha^-, \alpha^+, \beta^-, \beta^+)$. Then
$$
\lim_{n \to +\infty} e^{l( \alpha^n)+l( \beta^n)-l( \beta^n \alpha^n)}=[\alpha^-, \beta^-, \alpha^+, \beta^+]
$$
where  $l(\gamma)$ denotes the length of the closed geodesic corresponding to $\gamma$ for any $\gamma \in \Gamma$.
\end{proposition}

\noindent The set of  couples $(\alpha^-, \alpha^+)$ of all hyperbolic fixed points   of $\Gamma$ being dense in $X(\infty)^2$, this shows that  $[  f_\infty (\xi_1),   f_\infty (\xi_2),   f_\infty (\xi_3),   f_\infty(\xi_4)]=[\xi_1, \xi_2, \xi_3, \xi_4]$  $\forall \xi_1, \xi_2, \xi_3, \xi_4 \in X(\infty)$.
 
 \vspace{6mm}
 {\bf Proof of Proposition \ref{crossratio}.}\\
Fix $x \in X$.  
For  $n\geq 0$, set $\gamma_n=\beta^n\alpha^n$,   and let $\gamma_n^-$, $\gamma_n^+$ be  its  repelling and attractive fixed points. 
Consider  two sequences of points $a_k \in \; ]\alpha^-, \alpha^+[$ and $b_k \in \, ]\beta^-, \beta^+[$  such that  
 $\displaystyle \lim_{k \to +\infty} a_k=\alpha^-$ and $\displaystyle \lim_{k\to +\infty} b_k=\beta^+$; moreover, we can choose a sequence $n_k \rightarrow \infty$ such that  $\displaystyle \lim_{k \to +\infty} \alpha^{n_k}   a_k=\alpha^+$ and   $\displaystyle \lim_{k \to +\infty} \beta^{-n_k}   b_k=\beta^-$.\\
Now, for each $k$, let $B_k$ be a compact ball centered at $x$  containing $a_k$ and $b_k$.
Notice that   $\gamma_n^-$ and $\gamma_n^+$ tend respectively to $\alpha^-$ and $\beta^+$, so the sequence of geodesics $ ]\gamma_n^-, \gamma_n^+[ $  tend to  $]\alpha^-, \beta^+[$: namely, for $k$ fixed,  the distance between   $]\gamma_n^-, \gamma_n^+[$  and   $]\alpha^-, \beta^+[$, restricted to the compact set $B_k$, tends to $0$ when $n \rightarrow \infty$.  
We can then choose $n_k$ large enough so  that  
$$d\Bigl(]\gamma_{n_k}^-, \gamma_{n_k}^+[ \, \cap  B_k \; , \; ]\alpha^-, \beta^+[ \, \cap B_k\Bigr) <1/k$$
Call $a'_k$, $b'_k$ the projections of  $a_k,b_k$ on $]\gamma_{n_k}^-, \gamma_{n_k}^+[$; so, the sequences $(a'_k)_k$ and $(b'_k)_k$ also  converge  respectively to  $\alpha^-$ and $\beta^+$, and the sequences  $(\alpha^{n_k} a_k')_k   $, $(\beta^{-n_k}   b_k')_k   $  respectively  to  $\alpha^+$ and $\beta^-$.

\noindent We then have:
\begin{equation} \label{cross}
  [\alpha^- \!, \beta^- \!, \alpha^+ \!, \beta^+]
=   \lim_{k \rightarrow \infty}  \frac{  e^{  d(a_k',\alpha^{n_k}   a_k' )+d(\beta^{-n_k}   b_k' , b_k' )}  }
                                                   { e^{ d(a_k',b_k' )+d(\beta^{-n_k}   b_k' ,\alpha^{n_k}   a_k' ) }     }
=   \lim_{k \rightarrow \infty}  \frac{  e^{  d(a_k,\alpha^{n_k}   a_k )+d(\beta^{-n_k}   b_k, b_k )}  }
                                                   { e^{ d(a_k',b_k' )+d(\beta^{-n_k}   b_k' ,\alpha^{n_k}   a_k' ) }     }  
\end{equation}
by definition of the cross-ratio.
Notice that  the numerator gives exacty $e^{l( \alpha^{n_k})+l( \beta^{n_k})}$, as the points $a_k$ and $b_k$ lie on the axes of $\alpha, \beta$ respectively. On the other hand, for $k \gg 0$  
\begin{equation} \label{equality}
  d(a_k' , b_k' )+d(\beta^{-n_k}   b_k'  \; ,\; \alpha^{n_k}   a_k' )= l(\gamma_{n_k}).
\end{equation}
Indeed, let $V_\beta (b_k)$ be the hyperplane orthogonal to the axis of $\beta$,  passing through $b_k$. 
When $k$ is large enough, the point $\alpha^{n_k}   a_k'  $  is close to $\alpha^+$, in particular  it belongs to the half space bounded by  $\beta^{-n_k}(V_\beta (b_k))$   which contains $\beta^+$; consequently, the point $\gamma_{n_k}    a_k'  = \beta^{n_k}\alpha^{n_k}    a_k'  $ belongs to the half-space $V_\beta (b_k)$  which contains $\beta^+$,  so $b_k' $ lies on the geodesic $(\gamma_{n_k}^-, \gamma_{n_k}^+)$ between  $a_k' $ and  $\gamma_{n_k}    a_k'  $.  
As   $d( b_k' , \gamma_{n_k}   a_k' ) = d(\beta^{-n_k}   b_k' \,, \, \alpha^{n_k}   a_k' )$, the equality (\ref{equality}) readily follows. Letting $k \to +\infty$ in (\ref{cross}) then achieves the proof.$\Box$

\section{Entropy rigidity}\label{Entropyrigidity}

This section is devoted to the proof of   Theorem \ref{teor2}. \\
The proof is through the method of barycenter, initiated by  Besson-Courtois-Gallot \cite{bcg}, \cite{bcg2}, and follows the lines of  \cite{courtois} (Theorem 1.6,  holding for   compact quotients). \linebreak 
The main difficulty in the  finite volume, non compact case is to show that the map produced by the barycenter method is proper: we  will recall in \S\ref{rigidityproof} the main steps of the construction, referring the reader to  \cite{courtois} for the estimates which are now  well established,  while we will  focus on the new estimates necessary to prove properness.
For this, we will need accurates estimates  of the Patterson-Sullivan measure of some subsets of $X(\infty)$, which we will describe   in the first subsection.

\subsection{On the Patterson measure of non uniform lattices}
 
The  Patterson-Sullivan measures  of $\Gamma$ are a  family of finite measures $\mu=(\mu_x)$, indexed by points of $X$ and with support included in the limit set $L\Gamma \subset X(\infty)$, satisfying the following properties   (cp. for instance \cite{sullivan}, \cite{roblin} for details about their  construction):

\begin{enumerate}

\item they are absolutely continuous w.r. to each other: for any $x, x' \in X$   
\vspace{-3mm}

\begin{equation} \label{conformal}
\frac{d\mu_{x'}}{d\mu_x} (\xi)= e^{-\delta_{\Gamma}  b_\xi (x', x)}
\end{equation}

\item they are $\Gamma$-equivariant: for every $\gamma \in \Gamma$ and every Borel set $A \subset X (\infty)$
\vspace{-4mm}

\begin{equation} \label{gammainvariance}
  \mu_{x} (\gamma^{-1}A) = \mu_{\gamma   x} (A) 
\end{equation}
  \end{enumerate}
  
 \vspace{-2mm} 
\noindent When $\Gamma$  is a lattice, we will use the decomposition of $X$ explained in \S\ref{subsectiondecomposition} to describe the local behavior of Patterson-Sullivan  measures of $\Gamma$ on  the limit set $\Lambda_\Gamma =X(\infty)$. 
\vspace{1mm}

   For    $x \in X$ and $\zeta \in X (\infty)$, we consider the point $x\zeta(R)$ at distance $R$ from $x$ on the geodesic ray $[x,\zeta[$, and  define the  ``spherical cap'' $V_\zeta (x,  R) \subset X (\infty)$  as the set of points  of $ X(\infty)$ whose projection on the geodesic ray $[x,\zeta[$  falls between $x\zeta(R)$ and $\zeta$. \linebreak
   The proposition below gives an uniform estimate, which will be  crucial in the sequel,  for the measure  $\mu_{x}(V_\zeta (x,  R))$ of ``small'' spherical caps, i.e. when $R \gg 0$.
  \vspace{1mm}
   
\noindent So, let  ${\cal D} =  {\cal K} \cup  {\cal C}_1 \cup  \cdots\cup   {\cal C}_l$ be a decomposition of the Dirichlet domain of $\Gamma$ centered at some fixed point $x$,  corresponding to the maximal,  bounded parabolic subgroups $P_1, ...,  P_l$ of $\Gamma$ with fixed points $\xi_1, ...,  \xi_l$  as described  in   \ref{subsectiondecomposition}. 
  If $x\zeta(R)$ projects to the thick part $\bar {\cal K}$ of $\bar  X$, then formulas (\ref{conformal}) and  (\ref{gammainvariance})
easily give the  uniform lower estimate:  
\vspace{-6mm}

\begin{equation}\label{eqthick}
\mu_{x} ( V_\zeta (x,  R)) \stackrel{c}{\succeq} e^{-\delta_\Gamma R}
\end{equation}

\noindent (where $c$ is a positive constant, depending on the minimal mass of a spherical cap at distance less than $D=diam({\cal K})$ from $x$).  
On the other hand, when $x\zeta(R)$ projects to the cuspidal part, we  have:
 
 \pagebreak
\begin{proposition}
\label{proppatterson}
 There exists a constant    $c >0$  satisfying the following property. \linebreak   
Let $\zeta \!\in\! X (\infty)$  and assume that the point $x\zeta(R)$  belongs  to $\gamma   {\cal C}_i$,  $R>0$.
Then:
 \vspace{-3mm}
 
%
\begin{equation}
\label{estimate}
 \mu_{x} ( V_\zeta (x,  R))  \stackrel{c}{\succeq} \hspace{1mm} e^{-\delta_{\Gamma} (R+r)}   \hspace{1mm}v_{P_i} (x,2r)
 \end{equation}

\noindent  where   $r= b_{\xi_i}(x, \gamma^{-1}   x\zeta(R))$. 
\end{proposition}

 This estimate    stems  from a series of technical lemmas, and might be deduced from  work developed in \cite{peigne} and \cite{sch} (notice however that, in \cite{peigne}, $\mu_{x}$ has no atomic part, and in   \cite{sch} the parabolic subgroups are assumed to satisfy an additional, strong regularity assumption).  Since the estimate   is of independent interest, we will report for completeness the  proof of Proposition  \ref{proppatterson}, in full generality, in the Appendix.

\subsection{Entropy rigidity : proof of Theorem \ref{teor2} }
\label{rigidityproof}

Let $\bar X = \Gamma \backslash X$,   fix a point $x_0 \in X$ and call for short $b_{\xi} (x) = b_{\xi} (x, x_0)$. \\
The function $b_\xi$ is strictly convex if  $K_X \leq -1<0$,  since for every point $y$ we have:
 \vspace{-3mm}  
     
\begin{equation}
\label{eqhessian}
Hess_y \; b_{\xi} \geq g_y - (db_\xi)_y \otimes (db_\xi)_y
\end{equation}

\noindent where $g$ denotes the metric tensor of $X$; moreover,  it is known
 that if the equality holds in (\ref{eqhessian}) at every point $y$ and for every direction $\xi$,  then the sectional curvature is constant,  and $X$ is isometric to the hyperbolic  space $\mathbb{H}^{n}$. The idea of the proof is  to show that the condition $\delta_{\Gamma} = n-1$ forces the equality in (\ref{eqhessian}).

 Recall that,  for every measure $\mu$ on $X (\infty)$ whose support is not reduced to one point,  we can consider its  {\em barycenter},  denoted $bar[\mu]$,  that is the unique point of minimum of the function
$y \mapsto {\cal B}_\mu (y) = \int_{X (\infty)} e^{b_\xi ( y)} d\mu (\xi)$
(notice that this is $C^2$ and strictly convex function,   as  $b_{\xi}(y)$ is). If $supp(\mu)$ is not a single point,  it is easy to see   that $\lim_{y \rightarrow \xi} {\cal B}_\mu (y) = +\infty$ for all $\xi \in X (\infty)$ cp.  \cite{courtois}.
  
 Consider now  the map $F: X \rightarrow X$ defined by 
 \vspace{-4mm}
 
$$
F(x) = bar \left[ e^{-  b_{\xi} (x)} \mu_x\right] = argmin \left[ y \mapsto  \int_{X (\infty)} e^{b_\xi (y, x)}d\mu_x (\xi) \right] 
$$
where $(\mu_x)_x$ is the family of Patterson-Sullivan measures associated with the lattice $\Gamma$.

\noindent In \cite{courtois} it is proved that the map $F$ satisfies the following properties:
\vspace{2mm}
        
{\em 

{\bf (a)} $F$ is equivariant with respect to the action of $\Gamma$,  i.e. $F(\gamma   x) = \gamma F(x)$;

{\bf (b)}  $F$ is $C^2$,  with Jacobian:
 \vspace{-2mm}       
        
\begin{equation}
\label{eqjacobian}
   | Jac_x   F | \leq \left( \frac{\delta_{\Gamma}+1}{n}  \right)^n   det^{-1} ( k_x ) 
\end{equation}
}      

\noindent  where  $k_x (u, v)$ is the bilinear form on $T_xX$ defined as

\begin{equation}
\label{eqk}
 k_x  (u, v) 
 = \frac{ \int_{X (\infty)} e^{b_\xi (F(x), x))}  \left[ (db_\xi)_{F(x)}^2  +   Hess_{F(x)} b_\xi \right] (u, v) \; d\mu_x (\xi)}
            { \int_{X (\infty)} e^{b_\xi (F(x), x))}  \; d\mu_x (\xi)}
\end{equation}
\normalsize      

\noindent Notice that the eigenvalues of $k_x$ are all greater or equal than 1,  by (\ref{eqhessian}).

\noindent Property {\bf (a)} stems from the equivariance (i) of the family of Patterson-Sullivan measures with respect to the action of $\Gamma$,  and from the cocycle formula for the Busemann function: $b_\xi(x_0, x) +b_\xi(x, y)=b_\xi(x_0, y)$.
Property {\bf (b)} comes from the fact that the Busemann function is $C^2$ on Hadamard manifolds,  and is proved by direct computation,  which does not use cocompactness. \\
By equivariance,   the map $F$ defines a quotient map $\bar F: \bar X \rightarrow \bar X$,  which is homotopic to the identity through the homotopy 
$$
      \bar F_t (x)= bar \left[  e^{-b_\xi (x)} \left( t \mu_x + (1-t) \lambda_x \right) \right]  mod \; \Gamma,  \;\; t \in [0, 1]
$$
where $\lambda_x$ is the visual measure from $x$ (with total mass equal to the volume of $S^{n-1}$).

\noindent Actually,  the map 
$  F_t  =  bar   \left[  e^{-b_\xi (x)} \left( t \mu_x  +  (1 -  t) \lambda_x \right) \right] $ defines a map between   the quotient manifolds,  as it is still $\Gamma$-equivariant; moreover,   we have $bar   \left[ e^{-b_\xi (x)} \lambda_ x \right]   =   x$  since,  for all $ v \in T_xX$:  
$$
 \left( d \mathcal B_{e^{-b_\xi (x)} \lambda_ x} \right)_x (v)
= \int_{X (\infty)} (d b_\xi)_x (v) e^{ b_\xi (x)}e^{ -b_\xi (x)} d\lambda_x (\xi)
= \int_{U_{x}X} g_x (u,  v) du = 0.
$$

\vspace{3mm}     
We will now prove that:

\begin{proposition}
\label{propproper}
The homotopy map $\bar F_t$ is proper.
\end{proposition}

\noindent Assuming for a moment Proposition \ref{propproper},  the proof of Theorem \ref{teor2} follows by the degree formula: since $\bar F$ is properly homotopic to the identity,  it has degree one,  so
\vspace{-5mm}

\small
\begin{eqnarray*}   {\rm vol}(\bar X) =   \left| \int_{\bar X} \bar F^{\ast} dv_g  \right| &\leq&  \int_{\bar X} | Jac_x \bar F | dv_g 
\\
&\leq&
 \left( \frac{\delta_{\Gamma}+1}{ n }  \right)^n        \int_{\bar X} det^{-1} (k_x) dv_g 
 \\
&\leq& 
\left( \frac{\delta_{\Gamma}+1}{\delta_{\Gamma}(\mathbb{H}^n) +1 }  \right)^n       {\rm vol}(\bar X)
\end{eqnarray*}
\normalsize

\noindent as $det(k_x)\geq 1$ everywhere. So,  if $\delta_{\Gamma}=\delta_{\Gamma}(\mathbb{H}^n)=n-1$,  we deduce that 
$det(k_x)= 1$ everywhere and $k=g$,  hence the equality in the equation (\ref{eqhessian}) holds for every $y=F(x)$ and $\xi$. Since $F$ is surjective,  this shows that $X$ has constant curvature $-1$.$\Box$
 \vspace{5mm}       

{\bf Proof of Proposition \ref{propproper}}.\\
Denote by $\bar z$ the projection of a point $z  \in  X$ to $\bar X$, and set
 $\delta=\delta_{\Gamma}$; recall that $\delta=n-1$,  but we will  use this property only  at the end of the proof. \\
Let $\mu_x^t= e^{-b_\xi (x)} \left( t \mu_x + (1-t) \lambda_x \right)$: we need to show that if $t_k \rightarrow t_0$ and if $\bar x_k \rightarrow \infty$
in $\bar X$,  then $\bar y_k=\bar F_{t_k} (\bar x_k) = \overline{bar [ \mu_{x_k}^{t_k}]}$ goes to infinity too. \\
Now, assume by contradiction  that the points $\bar y_k$  stay in a compact subset of $\bar X$: so (up to a subsequence) $\bar x_k,  \bar y_k$ lift to points $x_k,  y_k$ such that   $y_k \rightarrow  y_0  \in X$ and $d(y_0,  x_k) = d(\bar y_0 ,  \bar x_k) = R_k \rightarrow \infty$.  \\ 
By the cocycle relation $b_{\xi} (y, x)  =  b_{\xi} (y, y_0)  +  b_{\xi}  (y_0,  x)$  and by  the density formula for the \nolinebreak Patterson-Sullivan measures $\frac{d\mu_x}{d\mu_{y_0}}(\xi) = e^{-\delta b_{\xi}(x, y_0)}$,  we have 
\vspace{-5mm}

\begin{eqnarray}
\label{eqintegrals}
    ( d{\cal B}_{\mu_x^t} )_y (v) 
    &=&  t  \int_{X (\infty)}  (db_{\xi})_y (v) \,    e^{b_\xi (y, y_0)}e^{(\delta+1) b_\xi (y_0, x)}  d\mu_{y_0}(\xi) \notag \\
\vspace{-2mm}  &\ & \qquad \qquad \qquad  \qquad + (1-t)    \int_{X (\infty)}   (db_{\xi})_y (v)    e^{ b_\xi (y, x)}  \,  d\lambda_x (\xi) 
\end{eqnarray}
We will now estimate the  two terms in  (\ref{eqintegrals})   and show that $( d{\cal B}_{\mu_{x_k}^{t_k}} )_{y_k}$ does not vanish for $R_k \gg 0$,  a contradiction.
So,  let   $\zeta_k$ be the endpoints of the geodesic rays $y_0x_k$   and  let 
$v_k=(\nabla b_{\zeta_k})_{y_k}$.  
Also,  consider the spherical caps $  V_{\zeta_k}(y_0,  R_k/2) $ and $ V_{\zeta_k}(y_0,  R_k)$. 

\noindent Let us first consider  the contributions of the two integrals of the right hand side  in (\ref{eqintegrals})   over $X \smallsetminus V_{\zeta_k}(y_0,  R_k/2)$.
 If $\xi \in X (\infty) \smallsetminus V_{\zeta_k}(y_0,  R_k/2)$,   the projection   of $\xi$ over $y_0\zeta_k$ falls closer to $y_0$ than to $x_k$,  hence $b_\xi (y_0, x_k) \leq  0$; moreover,  
$\vert b_\xi (y_k,  y_0)\vert  \leq d (y_k,  y_0) \rightarrow 0$,  so the first integral on $X \smallsetminus V_{\zeta_k}(y_0,  R_k/2)$ for   $x=x_k$,  $y=y_k$ and $v=v_k$  gives:
\vspace{-3mm}

$$
      \left|  \int_{ X \smallsetminus V_{\zeta_k}(y_0,  \frac{R_k}{2})}  (db_{\xi})_{y_k} (v_k)  e^{b_\xi (y_k, y_0)} e^{(\delta+1) b_\xi (y_0, x_k)}  \,  d\mu_{y_0}   \right| \leq  2  \parallel   \mu_{y_0} \parallel
      $$
for $k \gg 0$. Analogously,  the second integral on $X \smallsetminus V_{\zeta_k}(y_0,  R_k/2)$ yields
$$ 
       \left| \int_{ X \smallsetminus V_{\zeta_k}(y_0,  \frac{R_k}{2})}   (db_{\xi})_{y_k} (v_k)  e^{b_\xi (y_k, x_k)}   d\lambda_{x_k}  \right| \leq    2 \,  {\rm vol}(\mathbb S^{n-1})
$$
for $k \gg 0$,  since $|b_\xi (y_k,  x_k) - b_\xi (y_0,  x_k)| \leq d(y_k,  y_0)$. So,  these contributions are bounded.

\noindent We now compute the contributions of the  integrals  over $V_{\zeta_k}(y_0,  R_k/2) \smallsetminus  V_{\zeta_k}(y_0,  R_k)$.
For all $\xi \in V_{\zeta_k}(y_0,  R_k/2) $ we have that  $(\nabla b_\xi )_{y_0}    (\nabla b_{\zeta_k})_{y_0}$ is close to $1$,  for $R_k \gg0$; 
moreover,  as 
$$\left| (\nabla b_\xi )_{y_k}         v_k -   (\nabla b_\xi )_{y_0}          (\nabla b_{\zeta_k})_{y_0}  \right| \stackrel{k\rightarrow \infty}{\longrightarrow} 0,$$  we deduce that  $(db_{\xi})_{y_k} (v_k) >\frac12$ on  $V_{\zeta_k}(y_0,  R_k/2) $  for $k \gg 0$,  hence these contributions are positive. 

Finally,  let us compute the contributions of these integrals on the   caps $V_{\zeta_k}(y_0,  R_k) $.  
For  $\xi \in V_{\zeta_k}(y_0,  R_k) $,  consider the ray $y_0\xi$ from $y_0$ to $\xi$,  and the projection $P(t)$ on the geodesic $y_0\zeta_k$ of   the point $\xi(t):= y_0\xi(t)$. We have,    by  (\ref{eqoppositetriangle}) 
$$b_\xi(y_0,  x_k) 
\geq   \lim_{t \rightarrow \infty} [d(y_0,  P(t)) + d(P(t),  \xi(t)) ] - [ d(\xi(t),  P(t)) + d(P(t), x_k) ] - \epsilon 0
\geq R_k -  \epsilon $$   
  with $\epsilon = \epsilon(\pi/2)$.   Therefore we deduce that,  for $k \gg0$,   we have
 \begin{eqnarray}
\label{eqint1}
 \int_{_{V_{\zeta_k}(y_0,  R_k))}}   \hspace{-12mm} (db_{\xi})_{y_k} (v_k)         e^{b_\xi (y_k, y_0)} e^{(\delta+1) b_\xi (y_0, x_k)} \,  d\mu_{y_0} 
&\geq&    \frac14  e^{(\delta+1)(R_k -   \epsilon )} \mu_{y_0} (V_{\zeta_k}(y_0,  R_k)) \\
\label{eqint2}
 \int_{_{V_{\zeta_k}(y_0,  R_k))}} (db_{\xi})_{y_k} (v_k)    e^{b_\xi (y_k, x_k)}  \,  d\lambda_{x_k} 
& \geq &  \frac14   e^{(R_k -  \epsilon )}   \,  {\rm vol}(\mathbb S^{n-1}).
\end{eqnarray}
 It is clear that   the right-hand side of (\ref{eqint2}) goes to infinity when $R_k \gg0$; we will now prove that the right-hand side  of (\ref{eqint1}) also diverges for $R_k \rightarrow \infty$. This will conclude the proof,  as it will show that $d{\cal B}_{ \mu_{x_k}^{t^k} }(v_k)$ does not vanish for $k\gg0$ (being a convex combination of two positively  diverging terms).\\
So,  let   ${\cal D} =  {\cal K} \cup  {\cal C}_1 \cup  \cdots\cup   {\cal C}_l$ be a decomposition of the Dirichlet domain of $\Gamma$ centered at $y_0$  as in (\ref{decompositionfundamentaldomain}),   corresponding  to maximal,  bounded parabolic subgroups $P_1, ...,  P_l$  with fixed points $\xi_1, ...,  \xi_l $.
We know that $\bar x_k$ belongs to some cusp of $\bar X$,  so $x_k \in \gamma    {\cal C}_i$   for some $\gamma$;  let then $r_k = b_{\xi_i} (y_0,  \gamma^{-1}  x_k) \leq R_k$.\\
If  $\gamma   \xi_i$ falls in $V_{\zeta_k}(y_0,  R_k)$ and $\delta\gg0$,   as $K_X \geq -b^2$ we can use  Propositions  \ref{proppatterson} and   \ref{propVp}  to deduce that   
 \begin{eqnarray*}
 e^{(\delta+1)R_k} \mu_{y_0} (V_{\zeta_k} (y_0,  R_k))
 & \succeq &  e^{(\delta+1)R_k} e^{- \delta(R_k +r_k)}  
 \sum_{n \geq 0 } v_{P_i} (2r) e^{-\delta n}\\
 &   \succeq & e^{R_k- \delta r_k }       v_{P_i} (2 r_k )  \\  
  & \succeq&{e^{R_k- \delta r_k }   \over   \mathcal A_{P_i} ( x_0,  y_0,  r_k ) }.    
 \end{eqnarray*}

\noindent Since $K_X \leq -1$,  we know that $\mathcal A_{P_i} ( x_0,  y_0,  r_k )  \preceq e^{-(n-1)r_k}$, so we obtain
$$
e^{(\delta+1)R_k} \mu_{y_0} (V_{\zeta_k} (y_0,  R_k)) \succeq e^{R_k+(n-1-\delta)r_k}.
$$
On the other hand,  when $\gamma   \xi_i \not\in V_{\zeta_k}(y_0,  R_k)$,    we have,  by Propositions  \ref{proppatterson}  and  \ref{propVp}:
  $$ 
  \hspace{-1mm}   e^{(\delta+1)R_k} \mu_{y_0} (V_{\zeta_k}(y_0,  R_k))
  \; \stackrel{}{\succeq} \;  e^{ R_k -\delta r_k }   v_P (y_0,  2r_k )  
   \succeq  \;   e^{ R_k + (n-1 -\delta) r_k }.   
   $$
 Both lower bounds  tend to $+\infty$ as $k\to +\infty$, since $R_k \to +\infty$ and  $\delta\leq n-1$;  thus,  the integral in  (\ref{eqint1})  diverges. This concludes the proof that the map $\bar F_t$ is proper.$\Box$

\section{Appendix}

We report here, for completeness, a proof of the estimate given in Proposition \ref{proppatterson}.

\noindent To prove Proposition \ref{proppatterson}, we will need a series of elementary lemmas, where some equalities hold up to  some constant: so we will use the symbol $f \stackrel{C}{\approx} g $ to mean that two quantities $f$ and $g$ differ of at most  $C$.
To avoid cumbersome notations,   we will give the same name to these constants  in all the lemmas, meaning that they all hold  for the  choice of a  suitable  constant $C$ large enough. All these constants will depend on  the upper bound of the sectional curvature $K_X \leq -1$  and, possibly, on other parameters of $\bar X =\Gamma \backslash X$ which we will specify case by case.
\vspace{3mm}

 Recall that a parabolic group $P$ of isometries  fixing $\xi \in X (\infty)$  is called {\em bounded} if it  acts cocompactly on $X (\infty) - \{ \xi \}$ (as well as on every horosphere $\partial H$ centered at $\xi$). If ${\cal D}(P,x)$ is the Dirichlet domain of $P$ centered at $x$,  the sets 
  ${\cal S}_x= {\cal D}(P,x) \cap \partial H_\xi (x)$ and  the trace at infinity   ${\cal S}_x (\infty) = \overline{{\cal D}(P,x)} \cap X(\infty)$ of  ${\cal D}(P,x)$ 
 are compact, fundamental domains for the action of $P$ on $\partial H_\xi (x)$ and on $X (\infty)$, respectively.
\vspace{1mm}

\pagebreak

\noindent The following Lemmas  can be found, for instance,  in \cite{sch} (Lemmes  2.6, 2.7 and 2.9):

 \begin{lemma}\label{lemmadistance} There exists  a constant  $C>0$ with the following property. \\
Let $x \in X$ and $ \zeta \in X(\infty)$ be fixed. Then,   for any $\xi \in V_\zeta (x,R)$ we have:
\vspace{-3mm}

 $$d(x\zeta(R),x\xi(R)) \leq C.$$
\end{lemma}

\begin{lemma}\label{lemmacalotte}  There exists a constant  $C>0$  with the following property. \\
Let $x \in X$ and $ \zeta \in X(\infty)$ be fixed. Then:

\noindent (i)   for any $x'$ such that $d(x,x')< C$ we have
\vspace{-3mm}

 $$  V_\zeta (x', R + C) \, \subset  \, V_\zeta  (x, R)\, \subset  \,V_\zeta  (x', R - C) $$

\noindent (ii)   for any $\xi \in V_\zeta (x,R+2C)$ we have 
\vspace{-3mm}

$$  V_{\xi}  (x, R + C) \, \subset  \, V_\zeta  (x, R)\, \subset  \,V_{\xi}  (x, R - C) $$

\noindent provided that $R >C$.
\end{lemma}

\begin{lemma}\label{lemmacalotte-busemann}  
Let $P$ a bounded parabolic subgroup of $X$ fixing $\xi$, and let $S_x(\infty)$ as above.
There exists a constant  $C>0$ (depending on the diameter of ${\cal S}_x$)  with the following properties: for any $p \in P$ 
\vspace{1mm}

\noindent (i) if $d(x,p   x)>2R$, with $R>C$, then  $\forall \eta \in  S_x(\infty)$ we have $p   \eta \in V_\xi (x, R -C)$ and
\vspace{-3mm}

$$b_{p   \eta} (x_R, p x_R ) \stackrel{C}{\approx} d(x,p   x) - 2R$$

\noindent (ii) if $d(x,p   x)\leq 2R$, then  $\forall \eta \in  S_x(\infty)$ we have $p   \eta \in X(\infty) \setminus V_\xi (x, R +C)$ and
\vspace{-3mm}

$$b_{p   \eta} (x_R, px_R) \leq C$$

\noindent where $x_R=x\xi(R)$ is the point at distance $R$ from $x$ on the geodesic ray $[x,\xi[$.
\end{lemma}

 \vspace{2mm}
 {\bf Proof of Proposition  \ref{proppatterson}}. \\
Let  ${\cal D} =  {\cal K} \cup  {\cal C}_1 \cup  \cdots\cup   {\cal C}_l$ be a decomposition of the Dirichlet domain of $\Gamma$ centered at $x$,  corresponding to the maximal,  bounded parabolic subgroups $P_1, ...,  P_l$ of $\Gamma$  with  parabolic fixed points $\xi_1, ...,  \xi_l$, and with  ${\cal C}_i = {\cal D} \cap  H_{\xi_i}$, as described  in   \ref{subsectiondecomposition}. Moreover, let ${\cal S}_i (\infty) = \overline{{\cal D} } \cap X(\infty)$ be the fundamental domains for the action of $P_i$ on  $X(\infty)\setminus \{\xi_i\}$,
and let  $z_t := x\zeta(t)$ and  $x_{i,t}:= x \xi_i(t)$.

\noindent  {\em We assume that  $z_R$ belongs to  $\gamma  H_{\xi_i}$; so,  call for short  $\xi=\xi_i$, $P=P_i$,  
 $x_{R}=x_{i,R}$, ${\cal S}(\infty)={\cal S}_i (\infty)$ and  set $r=b_\xi (x, \gamma^{-1} z_R)$ hereafter.}
 
\noindent Now, first notice   that $  | b_\eta(x,  z_R)-R |$ is bounded, uniformly in  $\eta \in V_\zeta (x,  R)$, since for $t \gg 0$ we have
 \vspace{-4mm}
 
$$b_\eta (x, z_R)  
 \stackrel{\epsilon}{\approx} \left( d(x,z_R) + d(z_R, x\eta(t))\right)   - d(x\eta(t),z_R)  =R
$$ 
 for $ \epsilon=\epsilon (\frac{\pi}{2})$ as in (\ref{eqoppositetriangle}).
Thus,   the density formula (\ref{conformal}) yields 
\vspace{-3mm}

$$\mu_{x}\bigl(V_\zeta (x,  R)
\bigr)  
\stackrel{c}{\asymp} 
e^{-\delta_{\Gamma}R}  
 \mu_{ z_R }\bigl(V_\zeta (x,  R)\bigr) $$
 for some constant $c>0$  only depending on the upper bound  of the curvature. \pagebreak

\noindent  It is thus sufficient to show that 
\begin{equation}\label{tesi}
\mu_{ z_R }\bigl(V_\zeta (x,  R)\bigr) \succeq    e^{-\delta_\Gamma r}   v_P (x,2r)
 \end{equation} 
For this, we will analyse two different cases:
\vspace{4mm}

{\bf Case 1:  $\zeta \in \Gamma \xi$.}
${}$\\
$a)$ Assume first  $\gamma=1$, so $\zeta=\xi$  and  $z_R = x_R   \in H_\xi$.\\
We have, by Lemma \ref{lemmacalotte-busemann}:
\vspace{-6mm}

\begin{eqnarray}
\label{eq1} \nonumber
\mu_{ x_R }\left(V_\xi (x,  R)\right)  
 &\geq&  \mu_{ x_R } (\{\xi\})  
     + \hspace{-4mm} \sum_{\stackrel{p\in P }{p  {\cal S}(\infty) \subset V_\xi (x,  R)}} 
        \hspace{-4mm}  \mu_{ x_R }(p  {\cal S}(\infty) ) \\
&\geq&
\hspace{-6mm}  \sum_{\stackrel{p\in P }{d(x, p  x)\geq 2R+C}}
\hspace{-6mm}  \mu_{ x_R }(p  {\cal S}(\infty) ).
\end{eqnarray}
From the equivariance and  the density formula  (\ref{conformal}), (\ref{gammainvariance}) for   the family  $\mu_x$ we get
\begin{equation}\label{eq2}
\mu_{   x_R }(p  {\cal S}(\infty) ) 
= \int_{ {\cal S} (\infty) }e^{-\delta_{\Gamma}b_{p\eta}(  x_R, p   x_R )}\mu_{x_R}(d\eta)\notag  
   \asymp   \mu_{x_R }( {\cal S}  (\infty)) e^{ \delta_{\Gamma}(2R- d(x, p  x))}
\end{equation}
because $b_{p\eta}(  x_R, p    x_R )  \approx d(x,px) - 2R$, by Lemma \ref{lemmacalotte-busemann}. \\ 
Now,  $\mu_{ x_R}( {\cal S}  (\infty)) \asymp e^{-\delta_\Gamma R} \mu_{ x}({\cal S}(\infty)) \asymp e^{-\delta_\Gamma R} $, 
since $X(\infty) = P {\cal S}  (\infty) \cup \{ \xi \}$   and  the mass of $\mu_{x_R}$ is not reduced to one atom, 
so  $\mu_{x}( {\cal S}  (\infty))>0$. \\
Therefore, from (\ref{eq1}) we deduce that, for $\Delta$ large enough, we have 
\begin{equation}\label{estimation1cas1}
\mu_{ x_R}\left(V_{\xi} (x,  R)\right) 
 \;\; \succeq  \;\;  e^{\delta_\Gamma R}    \hspace{-15mm}
\sum_{\stackrel{p\in P_i}{  2R+C + \Delta \geq d(x, p  x)\geq 2R+C }} \hspace{-12mm}
e^{ -\delta_{\Gamma} d(x, p  x)} 
\;\; \succeq \;\; e^{-\delta_\Gamma R}   v_P (x,2R)
\end{equation} 
 as $v^\Delta_P (x,2R) \asymp v_P (x,2R)$ by Proposition \ref{propVp}, if $\Delta \geq \Delta_0$.
The estimate     (\ref{tesi})  follows in this case, since  $\zeta=\xi$ and $z_R=x_R$, so   $r=b_{\xi}(x,  z_R) = b_{\xi}(x,  x_R) = R$.
\vspace{1mm}

\noindent $b)$ Assume now that $\zeta =  \gamma    \xi$ for some   $\gamma \neq 1$.\\
We then set  $\xi'= \gamma   \xi = \zeta$, $x'=\gamma x$, $H_{\xi'} = \gamma H_\xi$,  $x'_{t}= \gamma x_{t}$ and $R'   := b_{\xi'} (x',z_R)$.\\
Notice that,   without    loss of generality, we can assume that $x'$ lies at  distance less than ${\rm diam}(\mathcal K)$ from  the geodesic ray  $x\xi'$ (actually, as $P$ acts cocompactly on $ \partial H_{\xi}$, we can replace $\gamma$ by $\gamma p$ for some suitable $p \in P$),   so $d(z_{R-R'}, x')$ and  $d(z_R, x'_{R'})$ are both bounded by $ 2{\rm diam}(\mathcal K) $.   \\
By Lemma  \ref{lemmacalotte}(i), there  exists   $C>0$ such that 
\vspace{-3mm}

$$ V_{\xi'} (x,R)  =  V_{\xi'} (z_{R-R'},R') \supset V_{\xi'} (x',  R'+C)$$ 
and then   (\ref{tesi})  follows from a), by applying the inequality (\ref{estimation1cas1})  to $\xi', x'$ and $R'$.
Actually, as $d(z_R, x'_{R'})$ is bounded, we have $d \mu_{x'_{R'}}/   d \mu_{z_R}   \asymp 1$ and we get from  (\ref{estimation1cas1})  
\vspace{-3mm}

$$
\mu_{z_R} ( V_{\xi'} (x, R) ) \succeq \mu_{x'_{R'}} ( V_{\xi'} (x', R'+C) )  
        \succeq       e^{-\delta_\Gamma R'} v_P (x', 2R')
$$
 
 \noindent which gives  (\ref{tesi}), since  in this case $\xi'=\zeta$ and  $r= b_\xi (x, \gamma^{-1} z_R) = b_{\xi'} (x', z_{R}) = R'$.
\vspace{4mm}

{\bf Case 2:  $\zeta \not\in \Gamma \xi$.} ${}$\\
$a)$ Assume first that $\gamma\!=\!1$, so  $z_R   \in H_\xi$.\\
Let  $C$ be the constant  in Lemma \ref{lemmacalotte}.  If    $ \xi \in V_\zeta (x,R  -2C )$  we call $ S=R-4C$, so  that $ \xi \in V_\zeta (x,S+2C)$ and we have   $V_\zeta (x,S) \supset V_\xi (x,S+C)$  by Lemma \ref{lemmacalotte}(ii). 
\linebreak
\noindent Notice that  we have  $d(z_{R}, x_{S}) \leq 5C$  by Lemma \ref{lemmadistance}.  
Therefore,  applying  again (\ref{estimation1cas1})  to $\xi, x$ and  $S$, we get
 \vspace{-3mm}
 
\begin{equation} \label{cista}
\mu_{ z_{R}}\bigl(V_\zeta (x,  {S})\bigr)   
    \succeq  \mu_{x_{S}} \bigl(V_{\xi}(x, S +C)\bigr)
    \succeq  e^{-\delta_\Gamma S}  v_P (x, 2S) 
 \end{equation}

\noindent and the estimate  (\ref{tesi})  follows, since here  $r=b_\xi (x, z_R) \approx b_\xi (x, x_{S}) =S $.
 \vspace{1mm}

\noindent  On the other hand, if $\xi \not\in V_\zeta (x,R -2C)$, let  $\bar \zeta$ be  the point at infinity of the geodesic  supporting $]\zeta, x]$, different from $\zeta$, and  let $\bar x$ be  the point of   $] \zeta, x] \cap \partial  H_{\xi}$ closest to  $\zeta$. \linebreak
Moreover, let  $\bar R:= d(\bar x,  z_R)$. Notice  that $z_R=\bar x \bar \zeta (\bar R)$  and that, 
setting  $\bar x_{\bar R} := \bar x \xi (\bar R)$,  
 we have $d( z_R, \bar x_{\bar R}) < C$, always by Lemma \ref{lemmadistance}, so $d \mu_{\bar x_{\bar R}}/   d \mu_{z_R}   \asymp 1$.\\
Now, we have $V_\zeta (x,  R)= X(\infty)\smallsetminus V_{\bar \zeta}(\bar x,  \bar R)$  and  
$X(\infty)\smallsetminus V_\zeta(x, R-2C) = V_{\bar \zeta}(\bar x,  \bar R+2C)$; \linebreak 
as $\xi \not\in V_\zeta (x,R -2C)$ we deduce that $V_\zeta (x,  R) \supset  X(\infty)\smallsetminus V_{\xi}(\bar x, \bar R-C )$ by  Lemma \ref{lemmacalotte}(ii).\\
 Hence,  
\vspace{-4mm}

\begin{equation}\label{uno}
\mu_{z_R}(V_\zeta (x,  R)) 
\succeq \mu_{\bar x_{\bar R} }\Bigl(X(\infty)\smallsetminus V_{\xi}(\bar x, \bar R-C )\Bigr).
\end{equation}
Similarly to case 1, we can estimate this  by applying  Lemma \ref{lemmacalotte-busemann} to $\bar x$ and $\xi$ :  
\vspace{-4mm}

\begin{eqnarray} \label{due}
\nonumber
\mu_{ \bar x_{\bar R}  }\Bigl(X(\infty) \smallsetminus V_{\xi} ( \bar x,  \bar R-C)\Bigr) 
&\geq&
\sum_{\stackrel{p\in P}{p {\cal S}(\infty) \, \cap  V_{\xi} (\bar x,  \bar R -C)=\emptyset}}
\mu_{ \bar x_{\bar R}  }(p  {\cal S} (\infty) )\\
\nonumber
&\geq& 
\sum_{\stackrel{p\in P}{d(\bar x, p  \bar x)\leq 2 (\bar R -2C)}}
\mu_{ \bar x_{\bar R}  }(p  {\cal S} (\infty) ) \\
&\asymp&
\sum_{\stackrel{p\in P}{d(\bar x, p  \bar x)\leq 2 (\bar R-2C)}}
\mu_{ \bar x_{\bar R}  }(  {\cal S} (\infty) ) 
\end{eqnarray}
 as $\mu_{ \bar x_{\bar R}  }(p  {\cal S} (\infty) ) \asymp \mu_{ \bar x_{\bar R}  }(  {\cal S} (\infty) ) $
 because, by Lemma \ref{lemmacalotte-busemann} $(ii)$,  
\begin{eqnarray*}
b_{p\eta}(  \bar x_{\bar R} , p    \bar x_{\bar R}  )&=&b_{p\eta}(  \bar x_{\bar R} ,     \bar x_{\bar R -2C}  )+b_{p\eta}(  \bar x_{\bar R -2C}, p   \bar  x_{\bar R -2C})+b_{p\eta}(  p    \bar x_{\bar R -2C}, p    \bar x_{\bar R}  )
\\
&\leq&4C+b_{p\eta}(  \bar x_{\bar R -2C}, p   \bar  x_{\bar R -2C})
\\  &\leq& 5C
\end{eqnarray*}
  if $d(\bar x ,p  \bar x) \leq 2 (\bar R -2C)$.  
 Moreover,  $d(z_R, \bar x_{\bar R})$ is bounded, so we deduce that 
$\mu_{\bar x_{\bar R}  }({\cal S} (\infty)) \asymp\mu_{ z_R  }(  {\cal S} (\infty) )  
  \asymp e^{-\delta_\Gamma  R}  $ and that
  \vspace{-4mm}
  
  $$\bar R  = b_\xi (\bar x, \bar x_{\bar R})  \approx b_\xi (\bar x, z_R)  = d( \partial H_\xi, \partial H_\xi(z_R)) \approx b_\xi (x, z_R) =r \;;$$
 so, combining (\ref{uno}) and (\ref{due})
 we obtain   
 \begin{equation} \label{noncista'} \mu_{z_R}(V_\zeta (x,  R))\succeq
e^{-\delta_\Gamma R} v_{P}(\bar x, 2\bar R ) \succeq    e^{-\delta_\Gamma r} v_{P}(x, 2 r )  
\end{equation}
(since $\bar x$ is at bounded distance from the orbit of $x$).
\vspace{1mm}

\noindent $b)$ Assume now that  $\gamma \neq 1$. \\
We  set   $\xi'= \gamma   \xi$, $x'=\gamma x$, $H_{\xi'} = \gamma H_\xi$,  $x'_{R}= \gamma x_{R}$, with $d( x', [x, \xi[ ) \leq diam ({\cal K})$,   and we proceed as above, according to the cases $\xi'  \in V_\zeta (x, R- 2C)$  or  $\xi'  \not\in V_\zeta (x, R- 2C)$. \\
 In the first case, we call $S :=R-4C$, $S' = b_{\xi'} (x',z_S)$, so  
 $V_\zeta (x, S) \supset V_{\xi'} (x,S+C)$ and we have $d(z_R, x'_{S'}) \leq 6 C + 2diam ({\cal K})$; then, using Lemma \ref{lemmacalotte}, we deduce  similarly to  (\ref{cista}),  that  \pagebreak
 \vspace{-4mm}
 
 $$\mu_{ z_{R}}\bigl(V_\zeta (x,  {S})\bigr)   
    \succeq  \mu_{x'_{S'}}  \bigl(V_{\xi'}(x', S' +2C)\bigr)
    \succeq  e^{-\delta_\Gamma S'}  v_P (x', 2S') 
$$
 
\noindent which yields  (\ref{tesi}), as here  $r=b_\xi (x, \gamma^{-1} z_R) \approx b_{\xi'} (x', z_S) =S' $.\\
In the second  case, we call again  $\bar \zeta$   the point at infinity opposite to $\zeta$ with respect to  $x$, \linebreak $\bar x$   the point of   $] \zeta, x] \cap \partial  H_{\xi'}$ closest to  $\zeta$, and we set 
$\bar R:= d(\bar x,  z_R)$,  $\bar x_{\bar R} := \bar x \xi' (\bar R)$. \linebreak
So,  $z_R=\bar x \bar \zeta (\bar R)$,  $d( z_R, \bar x_{\bar R}) < 2C$  and  $d \mu_{\bar x_{\bar R}}/   d \mu_{z_R}   \asymp 1$. Then,  we deduce as before that $V_\zeta (x,  R) \supset  X(\infty) \smallsetminus V_{\xi'}(\bar x, \bar R+C )$
and  we obtain,  analogously to (\ref{noncista'}), that 
$$ \mu_{z_R}(V_\zeta (x,  R))\succeq
e^{-\delta_\Gamma R} v_{P}(\bar x, 2\bar R )   $$
which concludes the proof as, in this case,   
 $\bar R  = b_{\xi'} (\bar x, \bar x_{\bar R})  \approx   b_{\xi'} (x', z_R) = r \;.\Box$



\end{document}